\documentclass[12pt]{article}

\title{On largest offsprings  \\
in a critical branching process \\
with finite variance} 
\author{Jean Bertoin\thanks{Institut f\"ur Mathematik, 
Universit\"at Z\"urich, 
Winterthurerstrasse 190, 
CH-8057 Z\"urich, Switzerland. \hfill \eject
Email: jean.bertoin@math.uzh.ch} }

\date{}

\usepackage{graphicx}
\usepackage{amsfonts,amsmath,amssymb,bbm}
\usepackage{setspace} 
\def\proof{\noindent{\bf Proof:}\hskip10pt}        
\def\QED{\hfill $\Box$}

\font\tenmath=msbm10 scaled 1200
\font\sevenmath=msbm7 scaled 1200
\font\Fivemath=msbm5 scaled 1200
\newfam\mathfam \textfont\mathfam=\tenmath
\scriptfont\mathfam=\sevenmath \scriptscriptfont\mathfam=\Fivemath

\def \\ { \cr }
\def\R{\mathbb{R}}
\def \1{1 \mkern -6mu 1} 
\def\N{\mathbb{N}}
\def\E{\mathbb{E}}
\def\P{\mathbb{P}}

\def\R{\mathbb{R}}

\def\Q{\mathbb{Q}}

\def \e{{\rm e}}

\def \d{{\rm d}}
\def \p{{\bf p}}

\newtheorem{theorem}{Theorem}

\newtheorem{lemma}{Lemma}

\begin{document}

\maketitle

\begin{abstract} 
We continue our study \cite{Be} of the distribution of the maximal number $X^{\ast}_k$ of offsprings amongst all individuals
in a critical Galton-Watson process started with $k$ ancestors, treating the case when the  reproduction law has a regularly varying tail $\bar F$
with index $-\alpha$ for $\alpha>2$ (and hence finite variance).
We show that  $X^{\ast}_k$ suitably normalized
converges in distribution to a Frechet law with shape parameter $\alpha/2$; this contrasts sharply with the case $1<\alpha <2$ when the variance is infinite. 
More generally, we obtain a weak limit theorem for the offspring sequence ranked in the decreasing order, in terms of atoms of a certain doubly stochastic Poisson measure. 

\end{abstract}

{\bf Key words:} Branching process, maximal offspring, extreme value theory, Cox process. 

{\bf Subject Classes:} 60F05, 60J80.
  
\begin{section}{Introduction and main results}
Let $\p=(p(n))_{n\in \N}$ be a probability measure on the set  of nonnegative integers; we view $\p$ as the law of a random variable $X$ which gives the number of children of a typical individual. We always assume criticality, i.e.
$$\E(X)=\sum_{n=0}^{\infty} n p(n) = 1\,,$$
and implicitly exclude the degenerate case $\p=\delta_1$. 
We  consider a Galton-Watson process having reproduction law $\p$. 
When the latter starts from $k\geq 1$ ancestors, we write $T_k$ for the size of the total population, and then for $i=1, \ldots, T_k$, we denote by $X_i$ the number of children of the $i$-th individual, where individuals in the branching process are enumerated by some arbitrary procedure (e.g. breadth first search). 
This note is concerned with weak limit theorems as $k\to \infty$ for the maximal offspring 
$$X^{\ast}_k=\max\{ X_i: 1\leq i \leq T_k\}\,,$$
and completes results of a previous work \cite{Be} which we now briefly recall. 
Plainly,  $X^{\ast}_k$ can be viewed as the maximum of $k$ i.i.d. copies of $X^{\ast}_1$, and one should naturally expect extreme value distributions to appear in the limit.

Let $\bar F(n)=\sum_{i=n+1}^{\infty}p(i)$ denote the tail distribution function of $\p$, and suppose that
\begin{equation}\label{E1}
\bar F(n) = n^{-\alpha} \ell(n)
\end{equation}
where $\ell$ is a slowly varying function and $\alpha\geq 1$. The main result in \cite{Be} is that for $1< \alpha<2$, $k^{-1}X^{\ast}_k$ converges in law as $k\to \infty$ to a Frechet distribution with shape parameter $1$. It is rather surprising that this weak limit theorem depends on  $\alpha$ only through the scale parameter of the Frechet law; in particular  in the normalization 
of $X^{\ast}_k$ is independent of $\alpha$ and the slowly varying function $\ell$ plays no role at all. 
The assumption $1<\alpha < 2$ implies that the variance
$$\sigma^2=\sum_{n=0}^{\infty} n^2p(n) -1$$
of the reproduction law is infinite, and it is therefore natural to wonder whether a similar phenomenon also occurs  in the case $\alpha >2$ for which the variance $\sigma^2$ is finite. The approach used in \cite{Be} for $1<\alpha < 2$ relies on the connexion  popularized by Harris between branching processes and left-continuous random walks, stable limit theorems for the latter, and the structure of jumps of stable L\'evy processes. It does not extend to the case $\alpha >2$ since Brownian motion  then arises in the limit and the jump structure becomes degenerate. 

In order to state the first result of the present contribution, it is convenient to introduce an asymptotic inverse of the tail distribution $\bar F$, that is a function $\varphi: (0,1)\to \R_+$ such that
$$\bar F(\varphi(\varepsilon))\sim \varepsilon\quad \hbox{as }\varepsilon\to 0+\,.$$
It is well-known that such a function exists when \eqref{E1} holds; further $\varphi$ then varies regularly at $0+$ with index $-1/\alpha$ and can be chosen nonincreasing. 

\begin{theorem}\label{T1} Assume \eqref{E1} with $\alpha>2$. Then for every $x>0$, we have
$$\lim_{k\to\infty} \P( X^{\ast}_k < x \varphi(k^{-2})) = \exp\left(-\frac{\sqrt 2}{\sigma}\, x^{-\alpha/2}\right)\,.$$
\end{theorem}

 Theorem \ref{T1} contrasts sharply with the results of \cite{Be}; this  points at the following transition for critical branching processes  with a reproduction law fulfilling \eqref{E1} and started from $k\gg 1$ ancestors. For $1<\alpha < 2$, the maximal offspring is always of order $k$, while for $\alpha>2$ it is only $k^{2/\alpha}$. More precisely, the maximal offsprings for two such branching processes
 may have much different asymptotic behaviors even when
 the two reproduction laws have the same finite variance, whereas  the asymptotic behaviors are always essentially the same when both variances are infinite.

It is well-known and easy to prove that when the reproduction law $\p$  has finite variance $\sigma^2$, the Laplace transform of $T_k$, the total population  generated by $k$ ancestors fulfills
$$\lim_{k\to \infty} \E(\exp(-ak^{-2}T_k))= \exp\left(-\sqrt{2a}/\sigma\right)\,,\qquad a\geq 0\,,$$
or equivalently that $k^{-2}T_k$ converges weakly towards some stable(1/2) variable $\tau$. Comparing with Theorem \ref{T1}, this suggest that it should be interesting to also take into account the role of $T_k$
in the study of maximal offsprings.

 In this direction, introduce  a sequence $(\eta_i)_{i\in\N}$ of i.i.d. variables with law $\p$, which we suppose further independent of the branching process, and hence of $T_k$.
 So $(X_i)_{i\in\N}$ and $(\eta_i)_{i\in\N}$ have the same law, but plainly  the partial sequences  $(X_i)_{i\leq T_k}$ and $(\eta_i)_{i\leq T_k}$ do not have the same distribution
since  $(X_i)_{i\in\N}$ and $T_k$ are not independent (and more precisely
there is the identity $T_k= X_1+\ldots + X_{T_k}+1$).  
Nonetheless,  it is easy to check that
$$\eta^{\ast}_{k}= \max \{\eta_i: 1\leq i \leq T_k\}\,.$$
 has the same asymptotical behavior as $X^{\ast}_k$ as $k\to \infty$, namely 
$$\lim_{k\to\infty} \P( \eta^{\ast}_k < x \varphi(k^{-2})) = 
\exp\left(-\frac{\sqrt 2}{\sigma}\, x^{-\alpha/2}\right)\,.$$
This reflects the fact that  when the reproduction law fulfills \eqref{E1} with $\alpha >2$, and the number of ancestors is large, the most prolific individual in a branching process has a negligible impact on the whole process. It is interesting to stress that this phenomenon ceases for $1<\alpha< 2$. More  precisely, in the case $1<\alpha < 2$, even though the maximal offspring $X^{\ast}_{k}$ is order $k$ and thus much smaller than the total population $T_k$ which has order 
$k^{\alpha}$,  $k^{-1} \eta^{\ast}_{k}$ and $k^{-1} X^{\ast}_{k}$ 
{\it do not} have the same asymptotic distribution; see Comment 2 in \cite{Be} . 

Our second result extends Theorem \ref{T1} by considering more generally not just the largest offspring, but also the second, third, ... largest, jointly with the total population size. Specifically, we write
$$X^{\ast}_{k,1}=X^{\ast}_k\geq X^{\ast}_{k,2}\geq X^{\ast}_{k,3}\geq \ldots$$
for the ordered sequence of the offspring numbers $\{X_i: 1\leq i \leq T_k\}$.

\begin{theorem}\label{T2} Assume \eqref{E1} with $\alpha>2$. Then 
$$\left(k^{-2} T_k; \frac{1}{\varphi(k^{-2})} X^{\ast}_{k,1}, \frac{1}{\varphi(k^{-2})} X^{\ast}_{k,2}\,,  \ldots , \,  \frac{1}{\varphi(k^{-2})} X^{\ast}_{k,T_k}\right)$$
converges in the sense of finite-dimensional distributions as $k\to \infty$ towards $(\tau;  \xi_1, \xi_2, \ldots)$, where
$\tau$ is a positive stable$(1/2)$ variable with law
$$\frac{\sigma}{\sqrt{2 \pi t^3}}\exp\left( - \frac{\sigma^2}{2t}\right)\d t\,, \qquad t>0\,, $$
and $\xi_1>\xi_2> \ldots$ denotes the sequence of the atoms ranked in the decreasing order of a random measure $M$ on $(0,\infty)$ such that conditionally on $\tau = t$, $M$ is Poisson with intensity $\alpha  t x^{-\alpha-1}\d x$. 
\end{theorem}
In particular, we see that the sequence of offspring numbers properly normalized and ranked in the decreasing order converges weakly towards the ranked sequence of the atoms of a doubly stochastic Poisson measure $M$ (also known as a Cox process). 

The rest of this work is organized as follows. Proofs of the two theorems are given in the next section, relying on estimates for Laplace transforms of linear functionals of offspring variables. We then gather in Section 3 a number of comments and observations, including some pointers at related literature.

\end{section}

\begin{section}{Proofs}

Our approach for establishing Theorems \ref{T1} and \ref{T2}  relies on standard techniques for estimating Laplace transforms of  linear functionals of the empirical measure of offsprings. One could also prove Theorem \ref{T2} by first rephrasing it in terms of random walks and then using more sophisticated tools of weak convergence of c\`adl\`ag semimartingales in Skorohod's space; however this would require knowing  {\it a priori}  the correct orders of magnitudes and also identifying   {\it a priori} dominant contributions and their relations. Therefore we prefer to develop arguments which are both more elementary and natural.

We also point out  that Theorem \ref{T1} can be recovered from Theorem \ref{T2} and the following elementary calculation. Taking  Theorem \ref{T2} for granted, we have for every $x>0$
\begin{eqnarray*}
\lim_{k\to\infty}\P(X^*_{k,1}\leq \varphi(k^{-2})x) &=&
\P(\xi_1\leq x)\\
& =& \P(M((x,\infty))=0)\\
&=& \E(\exp(- x^{-\alpha} \tau)) \\
&=& \exp\left(-\frac{\sqrt 2}{\sigma}\, x^{-\alpha/2}\right)\,.
\end{eqnarray*}
However it would probably not have been easy to guess {\it a priori} the correct renormalization
of the largest offsprings without knowing beforehand Theorem \ref{T1}.

Our starting point is based on the following elementary feature.
Consider a function $f:\N\to[0,\infty]$; we set
$$H(f)=\sum_{i=1}^{T_1}f(X_i).$$
and aim at evaluating  the Laplace transform 
$${\mathcal L}(f)=\E\left( \exp(-H(f))\right)\,.$$
In this direction, we observe from the branching property that ${\mathcal L}(f)$ solves the equation
\begin{equation}\label{E2}
s = \sum_{n=0}^{\infty} p(n) \e^{-f(n)} s^n\,, \qquad s\in[0,1].
\end{equation}
Since the entire function 
$$g_f: s\mapsto \sum_{n=0}^{\infty} p(n) \e^{-f(n)} s^n$$ 
has derivative $g'_f\leq 1$  on $[0,1]$ (recall that $\p$ is critical) and $0< g_f(0)< g_f(1)<1$, \eqref{E2} has a unique solution and thus determines ${\mathcal L}(f)$. 

\subsection{Proof of Theorem \ref{T1}}

The heart of the proof of Theorem \ref{T1} lies in the following lemma. At this stage, we do not require \eqref{E1} to hold.

\begin{lemma} \label{L1}  Let $\p$ be a critical reproduction law with finite variance $\sigma^2$. Suppose  that its  tail distribution has  $\bar F(x)>0$ for all $x>0$ and  that
$$\lim_{y\to \infty}\limsup_{x\to \infty} \frac{\bar F(xy)}{\bar F(x)}=0\,.$$
Then 
$$\P(X^{\ast}_1>x) \sim \sqrt{2\bar F(x)/Ê\sigma^2}\quad \hbox{as }x\to \infty.$$
\end{lemma}
{\bf Remark.} Note that $X^{\ast}_1$ may have infinite expectation; this occurs for instance when
$\bar F(x) \sim (x\ln x)^{-2}$ as $x\to \infty$. 

\proof Let $x>0$ and define $f(n)=0$ for $n\leq x$ and $f(n)=\infty$ for $n>x$. So 
$$H(f)=\left\{\begin{matrix}0\  \hbox{ if } X^{\ast}_1\leq x\\
\infty\  \hbox{ if } X^{\ast}_1> x \end{matrix} \right. $$
 and if we set ${\mathcal L}(f):= \varrho(x)$, then $\varrho(x)= \P(X^{\ast}_1\leq x)$ is the distribution function of $X^{\ast}_1$. In particular, \eqref{E2} shows that $\varrho(x)$ is the unique solution to the equation
\begin{equation}\label{E3}
\varrho(x) = \sum_{n\leq x} p(n) \varrho(x)^n\,,\qquad \varrho(x)\in[0,1]\,.
\end{equation}

Now introduce the generating function $g$ of $\p$, 
$$g(s)=\sum_{n=0}^{\infty} p(n) s^n\,, \qquad s\in[0,1]\,.$$
Since $\p$ is critical with variance $\sigma^2$, we have
$$g(s) = s + \frac{1}{2}\sigma^2 (1-s)^2 + o((1-s)^2)\,, \qquad s\to 1-\,.$$
Combining with \eqref{E3}, we obtain the estimate as $x\to \infty$
\begin{equation}\label{E4}
\frac{1}{2}\sigma^2 \bar \varrho(x)^2 + o(\bar \varrho(x)^2) = \sum_{n>x}p(n) \varrho(x)^n
= \bar F(x) \sum_{n=1}^{\infty} \frac{p(x+n)}{\bar F(x)} \varrho(x)^{n+x}\,,
\end{equation}
where $\bar \varrho=1-\varrho$ is the tail distribution of $X^{\ast}_1$. 
 We now claim that 
$$\sum_{n=1}^{\infty} \frac{p(x+n)}{\bar F(x)} \varrho(x)^{n+x}\sim 1\qquad \hbox{as }x\to \infty,$$
which will complete the proof of the statement, thanks to \eqref{E4}. 

Indeed, on the one hand, since  $({p(x+n)}/{\bar F(x)}, n\geq 1)$ is a probability measure and $\varrho(x)\leq 1$, there is the obvious upper bound
\begin{equation}\label{E*}
\sum_{n=1}^{\infty} \frac{p(x+n)}{\bar F(x)}\varrho(x)^{n+x} \leq 1\,.
\end{equation}
On the other hand, for every fixed $y>0$, we have
\begin{equation}\label{E5}
\sum_{n=1}^{yx} \frac{p(x+n)}{\bar F(x)} \varrho(x)^{n+x}\geq \frac{\bar F(x)-\bar F(x(1+y))}{\bar F(x)}Ê\, \varrho(x)^{x(1+y)}\,.
\end{equation}
We now see from \eqref{E4} and \eqref{E*} that 
$$
 \bar \varrho(x)  \leq \sqrt{2\bar F(x)/\sigma^2} + o( \bar \varrho(x))\,.
$$
Further, $\bar F(x) = o(x^{-2})$ because $\p$ has a finite variance. Hence $\bar\varrho(x)=o(1/x)$, which in turn entails that 
$$
\lim_{x\to \infty}\varrho(x)^{(1+y)x} = 1\,.
$$
Plugging this in  \eqref{E5}, we get  that
$$\liminf_{x\to \infty} \sum_{n=1}^{\infty} \frac{p(x+n)}{\bar F(x)} \varrho(x)^{n+x} \geq 
\liminf_{x\to \infty} \frac{\bar F(x)-\bar F(x(1+y))}{\bar F(x)}\,,$$
and since $y$ can be chosen arbitrarily large, we conclude from  the assumption of the lemma that
$$\liminf_{x\to \infty} \sum_{n=1}^{\infty} \frac{p(x+n)}{\bar F(x)} \varrho^{n+x}(x) \geq 
1\,,$$
which ends the proof.
\QED

If we now assume that \eqref{E1} holds with $\alpha>2$, then the conditions of Lemma \ref{L1} are fulfilled. Theorem \ref{T1} then follows from the fact that, because different ancestors produce i.i.d. branching processes, $X^{\ast}_k$ can be viewed as the maximum of $k$ i.i.d. copies of $X^{\ast}_1$, the 
estimate of Lemma \ref{L1} and the classical result of Gnedenko (see, for instance Proposition 1.11 in \cite{Resnick}).

\subsection{Proof of Theorem \ref{T2}}

We consider a continuous function $f:\R_+\to \R_+$
with $f\equiv 0$ on some neighborhood of $0$, and $a\geq 0$. Recall that $\varphi$ is an asymptotic inverse of $\bar F$. 
For every $k\geq 1$, we write 
$$f_{k,a}(x)=ak^{-2} + f( x/\varphi(k^{-2}))\,.$$ The key step consists in estimating 
$$\bar {\mathcal L}(f_{k,a})= 1- {\mathcal L}(f_{k,a}) = 1- \E\left( \exp\left(-a k^{-2}T_1 - \sum_{i=1}^{T_1} f( X_i/\varphi(k^{-2}))\right) \right)$$
as $k\to \infty$. 

\begin{lemma}\label{L2} Assume \eqref{E1} with $\alpha>2$. 
In the notation above, we have
$$\lim_{k\to \infty} k\bar {\mathcal L}(f_{k,a}) =  \frac{\sqrt 2  }{\sigma}\left( a+ \alpha \int_0^{\infty} \left(1-\e^{-f(x)} \right) x^{-\alpha-1} \d x\ \right)^{1/2} \,, \qquad k \to \infty.$$
\end{lemma}

\proof For the sake of simplicity, we assume that $f\equiv 0$ on $[0,1]$, the general case only requiring slightly heavier notation. 
The calculations are closely related to those in the proof of Lemma \ref{L1}, and we shall therefore sometime provide fewer details. We first recall from \eqref{E2} that there is the identity
$$ {\mathcal L}(f_{k,a}) = \sum_{n=0}^{\infty} p(n) \exp(-ak^{-2}-f( n/\varphi(k^{-2})))  {\mathcal L}(f_{k,a})^n\,,$$
and then, using the expansion of the generating function $g$, we get that
\begin{equation}\label{E6}
\frac{1}{2}\sigma^2 \bar  {\mathcal L}(f_{k,a})^2 + o(\bar  {\mathcal L}(f_{k,a})^2) = \sum_{n=0}^{\infty}
p(n) \left(1-  \exp(-ak^{-2}-f( n/\varphi(k^{-2}))) \right)   {\mathcal L}(f_{k,a})^n \,.
\end{equation}

We first use that $f\equiv 0$ on $[0,1]$ and $ {\mathcal L}(f_{k,a}) \leq 1$, and get the upper bounds
\begin{eqnarray*}
\frac{1}{2}\sigma^2 \bar  {\mathcal L}(f_{k,a})^2 + o(\bar  {\mathcal L}(f_{k,a})^2) 
&\leq &  \sum_{n\leq \varphi(k^{-2})}
p(n)\left(1-  \e^{-ak^{-2}} \right) +  \sum_{n>\varphi(k^{-2})} p(n)\\
&\leq & a k^{-2} + \bar F(\varphi(k^{-2})) \,.
\end{eqnarray*}
Since $\bar F(\varphi(k^{-2})) \sim k^{-2}$, we conclude that $\bar  {\mathcal L}(f_{k,a})= O(1/k)$.
Recall that $\varphi$ is regularly varying at $0+$ with index $-1/\alpha >-1/2$, in particular $\varphi(k^{-2})= o(k)$ and  therefore
\begin{equation}\label{E7}
\lim_{k\to \infty}  {\mathcal L}(f_{k,a})^{\varphi(k^{-2})}=1\,.
\end{equation}

We need to estimate the right hand side of \eqref{E6} as $k\to \infty$; in this direction it is convenient to decompose the sum depending on whether  $n\leq \varphi(k^{-2})$ or $n> \varphi(k^{-2})$,  as the summand has different asymptotic behaviors on these two regions. So first, as $f\equiv 0$ on $[0,1]$, we deduce from \eqref{E7} that
\begin{eqnarray*}
 &&\sum_{n\leq \varphi(k^{-2})}
p(n) \left(1-  \exp(-ak^{-2}-f( n/\varphi(k^{-2}))) \right)   {\mathcal L}(f_{k,a})^n \\
&\sim&   \sum_{n\leq \varphi(k^{-2})}
p(n) \left(1-\e^{-ak^{-2}}\right)
\sim   a k^{-2} \,.
\end{eqnarray*}

Second,  the assumption \eqref{E1}  ensures  that for every bounded and continuous function $h:[1,\infty)\to \R$, one has
$$\sum_{n> \varphi(k^{-2})}
p(n)  h(n/\varphi(k^{-2}) )  \sim  \bar F( \varphi(k^{-2})) \alpha
 \int_1^{\infty} h(x) x^{-\alpha-1} \d x\,.$$
Since ${\mathcal L}(f_{k,a})\leq 1$ and $\bar F( \varphi(k^{-2})) \sim k^{-2}$, we then readily deduce from \eqref{E7} that
 $$\sum_{n> \varphi(k^{-2})}
p(n) \left(1-  \exp(-ak^{-2}-f( n/\varphi(k^{-2}))) \right)   {\mathcal L}(f_{k,a})^n \sim
 k^{-2} \alpha 
 \int_0^{\infty} \left(1-\e^{-f(x)} \right) x^{-\alpha-1} \d x\,.$$
These two estimates combined with \eqref{E6} complete the proof of the lemma. \QED 

We can now tackle the proof of Theorem \ref{T2}. Because the branching processes generated by distinct ancestors are independent, we have
$$\E\left(\exp\left\{-ak^{-2}T_k - \sum_{n=1}^{T_k}f(X_i/\varphi(k^{-2}))\right\}\right)
= {\mathcal L}(f_{k,a})^k\,,$$
and it follows from Lemma \ref{L2} that 
\begin{eqnarray} \label{E8}
& &\lim_{k\to \infty}
\E\left(\exp\left\{-ak^{-2}T_k - \sum_{n=1}^{T_k}f(X_i/\varphi(k^{-2}))\right\}\right) \nonumber \\
&=&
\exp\left\{- \frac{\sqrt 2  }{\sigma}\left( a+ \alpha \int_0^{\infty} \left(1-\e^{-f(x)} \right) x^{-\alpha-1} \d x\ \right)^{1/2}
\right\}\,.
\end{eqnarray}

We next identify the right hand side of \eqref{E8}  as the joint Laplace transform of $(\tau, \langle M, f \rangle )$ where $\tau$ is a random variable and $M$ a random point measure on $(0,\infty]$ whose joint  law is specified in Theorem \ref{T2}. 
Indeed, we have then, first,
$$\E(\exp(-a \tau))= \exp(-\sqrt{2a/\sigma^2})\,,$$
second, 
$$\E(\e^{-\langle M,f\rangle}\mid \tau=t)=\exp\left\{-t \alpha \int_0^{\infty} (1-\e^{-f(x)}) x^{-\alpha-1} \d x\right\}\,,$$
and finally
$$\E(\exp(-a \tau-\langle M,f\rangle))= \exp\left\{- \frac{\sqrt 2  }{\sigma}\left( a+ \alpha \int_0^{\infty} \left(1-\e^{-f(x)} \right) x^{-\alpha-1} \d x\ \right)^{1/2}
\right\}\,.$$

Next introduce for each $k\geq 1$ the empirical measure of rescaled offsprings when the branching process has $k$ ancestors
$$M'_k=\sum_{i=1}^{T_k} \delta_{X_i/\varphi(k^{-2})}\,.$$
We ignore the atom of $M'_k$  at $0$ to get a random point measure $M_k$ on $(0,\infty]$, and  rephrase \eqref{E8} as
$$\lim_{k\to \infty} \E\left(\exp\left\{-a k^{-2}T_k - \langle M_k, f\rangle\right\}\right) 
= \E(\exp(-a \tau-\langle M,f\rangle))$$
for all $a\geq 0$ and continuous functions $f:(0,\infty]\to \R_+$ which have compact support. According to Theorem 14.16 in \cite{Kallenberg}, this ensures that the pair $(k^{-2}T_k, M_k)$ converges weakly as $k\to \infty$ towards $(\tau, M)$, where the space of Radon measures on $(0,\infty]$ is endowed with the vague topology. 
This entails Theorem \ref{T2}. 

\end{section}

\begin{section}{Miscellaneous remarks}

1. 
The literature surveyed by  Yanev \cite{Yanev1, Yanev2} contains a variety of results about extremes for individuals in branching processes. In particular, Rahinov and Yanev \cite{RY} have characterized the asymptotic behavior as $n\to\infty$ of the maximal offspring at the $n$-th generation conditionally on the event that the extinction has not yet occurred; previously Pakes \cite{Pakes} considered asymptotics of the largest score up to and including the $n$-th generation, again conditionally on the event that extinction does not occur before the $n$-th generation. 
In a different but related direction, namely the distribution of the maximum of the branching process itself, we refer to \cite{Ath,  BV, Lind, VWF}.

2. If one replaces  the assumption \eqref{E1} in  Theorem \ref{T1} by the condition that the tail distribution of the reproduction law decays exponentially, say $\bar F(x)\sim a \e^{-bx}$ for some $a,b>0$, then the same argument shows that 
$X^{\ast}_k-2b^{-1}\ln k$ now converges weakly as $k\to \infty$ to a Gumbel distribution. 
In the same way, the analogue of Theorem \ref{T2} can be stated as follows: As $k\to \infty$,  
$$\left(k^{-2} T_k;  X^{\ast}_{k,1}-2b^{-1}\ln k, X^{\ast}_{k,2}-2b^{-1}\ln k\,,  \ldots , \,   X^{\ast}_{k,T_k}-2b^{-1}\ln k\right)$$
converges in the sense of finite-dimensional distribution towards $(\tau;  x_1, x_2, \ldots)$, where
$\tau$ has the law
$$\frac{\sigma}{\sqrt{2 \pi t^3}}\exp\left( - \frac{\sigma^2}{2t}\right)\d t\,, \qquad t>0\,, $$
and conditionally on $\tau = t$, 
$x_1>x_2> \ldots$ denotes the sequence of the atoms ranked in the decreasing order of a Poisson measure on $(-\infty,\infty)$ with intensity $t ab \e^{-bx }\d x$. 

3. We point out that Equation \eqref{E2} has a natural interpretation in terms of branching processes. In this direction, we may view $\p_f=(p_f(n)=p(n)\e^{-f(n)})_{n\in\N}$ as a probability measure on $\N\cup\{\infty\}$ by setting $p_f(\infty)=1-\sum_{n=0}^{\infty}p(n)\e^{-f(n)}$. Then 
$$\sum_{n=0}^{\infty} p(n) \e^{-f(n)} s^n\,, \qquad 0\leq s < 1$$ 
should be thought of as the generating function of $\p_f$, and if we consider a Galton-Watson branching process with values in $\N\cup\{\infty\}$ and reproduction law $\p_f$ (of course this process is absorbed at $\infty$ as soon as an individual has an infinite offspring), then we recognize  \eqref{E2} as the equation for the probability of extinction. 

In particular, the well-known formula of Dwass \cite{Dwass} for the distribution of the total population size yields the semi-explicit formula
 $${\mathcal L}(f)=\sum_{n=1}^{\infty} \frac{1}{n} \p^{\ast n}_f(n-1)\,,$$
where  $ \p_f^{\ast j}= \p_f\ast \ldots \ast \p_f$ denotes the $j$-th convolution power of $\p_f$. Unfortunately this expression seems too complicated to be of some use in practice. 

4. It is natural to view $H(f)$ as a Hamiltonian, ${\mathcal L}(f)$ as a partition function, and
 introduce the Gibbs measure 
$$\Q=\frac{\e^{-H(f)}}{{\mathcal L}(f)}\P\,,$$
where $\P$ denotes the law of the Galton-Watson process with reproduction law $\p$ and started from a single ancestor. Using \eqref{E2} and for instance Otter's formula (see, e.g. Formula (304) in Pitman \cite{PiSF}), it is easily seen that $\Q$ is again the law of a Galton-Watson process with reproduction law given by
$$q(n)= {\mathcal L}(f)^{n-1} \e^{-f(n)} p(n)\,,\qquad n\in \N\,.$$
In particular, it is now elementary to compute the relative entropy $D(\Q\|\P)$ of $\Q$ with respect to $\P$. Indeed, excluding implicitly the trivial case $f\equiv 0$, Wald's lemma yields 
$$\E_{\Q}(H(f)) =  \langle {\bf q}, f\rangle \E_{\Q}(T_1) \,,$$
with $ \langle {\bf q}, f\rangle=Ê\sum_{n=0}^{\infty}f(n) q(n)$, and
it is a standard fact for sub-critical branching processes that 
$$1/\E_{\Q}(T_1)=1-\E_{\Q}(X_1)= 1-\sum_{n=0}^{\infty} n q(n)\,.$$

5. It would be interesting to obtain an analogue of Theorem \ref{T2} in the case $1<\alpha < 2$. As  a first modest step, one should probably try to understand the impact of the maximal offspring $X^{\ast}_k$ on
the total population $T_k$.

\end{section}

  \end{document}